\documentclass[11pt]{amsart}
\input epsf
\usepackage{latexsym}
\usepackage{amssymb,amsmath,amsthm,amscd, amssymb,
graphicx, enumerate, verbatim}
\usepackage[all]{xy}

\numberwithin{equation}{section}
\newtheorem{cor}[equation]{Corollary}

\newtheorem{thm}[equation]{Theorem}

\newtheorem{Example}[equation]{Example}

\newtheorem{remark}[equation]{Remark}
\newenvironment{rmk}{\begin{remark}\rm}{\end{remark}}
\def\co{\colon\thinspace}

\newcommand{\diam}{\mbox{diam}}

\newcommand{\Iso}{\mbox{Iso}}

\newcommand{\pinch}{\mbox{pinch}}
\newcommand{\cd}{\mbox{cd}}
\newcommand{\e}{\epsilon}
\def\a{\alpha}
\def\G{\Gamma}

\def\d{\delta}
\def\r{\rho}
\def\s{\sigma}
\def\l{\lambda}

\newcommand{\g}{\gamma}

\begin{document}

\abovedisplayskip=6pt plus3pt minus3pt
\belowdisplayskip=6pt plus3pt minus3pt

\title[Pinching estimates for nilpotent groups]{\bf
Pinching estimates for negatively curved manifolds with nilpotent
fundamental groups}

\rm

\author{Igor Belegradek}
\address{Igor Belegradek\\School of Mathematics\\ Georgia Institute of
Technology\\ Atlanta, GA 30332-0160}\email{ib@math.gatech.edu}
\author{Vitali Kapovitch}
\address{Vitali Kapovitch\\Department of Mathematics\\University of California\\
Santa Barbara, CA 93110}\email{vitali@math.ucsb.edu}
\date{}
\begin{abstract}
Let $M$ be a complete Riemannian metric of sectional curvature
within $[-a^2,-1]$ whose fundamental group contains a
$k$-step nilpotent subgroup of finite index. We prove that $a\ge k$
answering a question of M.~Gromov.
Furthermore, we show that for any $\epsilon>0$,
the manifold $M$ admits a complete Riemannian metric of sectional curvature
within $[-(k+\epsilon)^2,-1]$.
\end{abstract}
\maketitle

\section{Introduction}

If the fundamental group of a complete pinched negatively curved
manifold is amenable, it must be finitely generated and virtually
nilpotent~\cite{BuSc, Bow, BGS}. In this paper we relate the
nilpotency degree of the group to the  pinching of the negatively
curved metric.

\begin{thm} \label{thm: bound on pinch}
Let $M$ be complete Riemannian manifold
with sectional curvature satisfying
$-a^2\le \sec(M)\le -1$. If $\G$ is a $k$-step nilpotent
subgroup of $\pi_1(M)$, then $a\ge k$. In particular,
if $a\in [1,2)$, then $\G$ is abelian.
\end{thm}

If the cohomological
dimension $\cd(\G)$ of $\G$ equals to
$\dim(M)-1$, which if $\dim(M)>2$ is equivalent to assuming that
$\G$ acts cocompactly on horospheres,
Theorem~\ref{thm: bound on pinch} follows from
the proof of Gromov's theorem of almost flat manifolds
(see~\cite[Corollary 1.5.2]{BuKa}), by combining the
commutator estimate in almost flat
horosphere quotients with the displacement estimate
coming from the exponential convergence of geodesics.

More recently,  Gromov sketched in~\cite[p.309]{Gro}
a proof of the more general estimate 
\[a\ge \frac{k}{r+1}\ \ \mathrm{for}\ \
r=\left[\frac{\dim(M)-1-\mathrm{cd}(G)}{2}\right],\]
where $[x]$ denotes the largest integer satisfying $\le x$.
If $k\le r+1$, the estimate gives no information, so
Gromov asked~\cite[p.309]{Gro} whether it can be improved to
an estimate that is nontrivial for all $\cd(G)<\dim(M)$.
Theorem~\ref{thm: bound on pinch} provides a satisfying
answer that involves no dimension assumptions whatsoever.
The proof of Theorem~\ref{thm: bound on pinch} follows
the original Gromov's idea in~\cite{BuKa}, except that the
commutator estimate is run in a ``central'' orbit of an
N-structure given by the collapsing theory of J.~Cheeger,
K.~Fukaya, and Gromov~\cite{CFG}.
In ~\cite{BK} we proved the following classification theorem:

\begin{thm}\label{intro: maincor}\cite{BK}
A smooth manifold $M$ with amenable fundamental group admits a
complete metric of pinched negative curvature if and only if it is
diffeomorphic to the M\"obius band, or to the product of a line
and the total space a flat Euclidean vector bundle over a compact
infranilmanifold.
\end{thm}

The "if" direction in Theorem~\ref{intro: maincor} involves an
explicit warped product construction of a negatively pinched
metric on the product of $\mathbb R$ and the total space of a
flat Euclidean bundle over a closed infranilmanifold. 
By improving this warped product construction, we 
show that the pinching bounds provided by
Theorem~\ref{thm: bound on pinch} are essentially optimal.

\begin{thm}\label{optimal warp}
If $M$ be a pinched negatively curved manifold
such that $\pi_1(M)$ has a $k$-step nilpotent
subgroup of finite index, then $M$ admits a
complete Riemannian metric
of $\sec(M)\in [-(k+\e)^2, -1]$ for any $\e>0$.
\end{thm}

The metric constructed in Theorem~\ref{optimal warp}
has cohomogeneity one, specifically $M/\Iso(M)$ is
diffeomorphic to $\mathbb R$ (with the only exception
when $M$ is the M\"obius band equipped with a hyperbolic
metric).

We do not know whether $M$ in Theorem~\ref{optimal warp} always
admits a complete metric with $\sec(M)\in [-k^2, -1]$. This does
happen for $k=1$, since as we show in~\cite{BK} any complete
pinched negatively curved manifolds with virtually abelian
fundamental group admits a complete hyperbolic metric.

Another way to phrase the optimality of Theorem~\ref{thm: bound on
pinch} is via the concept of pinching. Given a smooth manifold
$M$, we define $\pinch^\mathrm{diff}(M)$ to be the infimum of
$a^{2}\ge 1$ such that $M$ admits a complete Riemannian metric of
$-a^2\le \sec(M)\le -1$. 
If $M$ admits no complete metric of
pinched negative curvature, it is convenient to let
$\pinch(M)^\mathrm{diff}=+\infty$. 
We then define
$\pinch^\mathrm{top}(M)$ to be the infimum of all
$\pinch^\mathrm{diff}(N)$ where $N$ is homeomorphic to $M$, and
define $\pinch^\mathrm{hom}(M)$ to be the infimum of
$\pinch^\mathrm{diff}(N)$'s where $N$ is manifold with
$\dim(N)=\dim(M)$ that is homotopy equivalent to $M$. Of course,
$\pinch^\mathrm{diff}(N)\ge \pinch^\mathrm{top}(M)\ge
\pinch^\mathrm{hom}(N)\ge 1$.

In general, the pinching invariants are hard to estimate and even
harder to compute (see~\cite{Gro}
and~\cite[Section 5]{Bel-inv} for surveys). 
Combining Theorems~\ref{thm: bound on
pinch}--\ref{optimal warp}, we compute the invariants in case
$\pi_1(M)$ is virtually nilpotent.

\begin{cor}
If $M$ be a pinched negatively curved manifold such that
$\pi_1(M)$ has a $k$-step nilpotent subgroup of finite index, then
$\pinch^\mathrm{diff}(M)=
\pinch^\mathrm{top}(M)=\pinch^\mathrm{hom}(M)=k^{2}$.
\end{cor}

This work was partially supported by the NSF grants
\# DMS-0352576 (Belegradek) and \# DMS-0204187 (Kapovitch).
We are thankful to J.~Cheeger for a discussion on collapsing.

\section{Proof of Theorem~\ref{thm: bound on pinch}}
\label{sec: proof}
A Riemannian metric is called $A$-{\it regular} if
$A=\{A_i\}$ is a sequence of nonnegative reals such that
the norm of the curvature tensor satisfies
$||\nabla^i R||\le A_i$. We call a metric
{\it regular} if it is $A$-regular for some $A$.
The collapsing theory works best for regular metrics,
and the Ricci flow can be used to deform any metric
with bounded sectional curvature to a
complete Riemannian metric that is close to the
original metric in uniform $C^1$ topology, has almost
the same sectional curvature bounds, and is regular.
(This fact has been known to some experts, but the first
written account only recently appeared in~\cite{Kap-ricflow}).
Thus we fix an arbitrary $\d>0$ and replace the given metric on $M$
by a nearby $A$-regular metric $g$ with $\sec_g\in [-(a+\d)^2,-1]$,
and then prove that $a+\d\ge k$, which would
imply $a\ge k$ because $\d$ is arbitrary.

Since the Riemannian covering of $(M,g)$
corresponding to $\G\le\pi_1(M)$
has the same curvature bounds as $(M,g)$, we can assume
that $\pi_1(M)=\G$. Denote the universal cover of $M$ by $X$.
If $k=1$, all we assert is $a\ge 1$ which is trivially true,
so we assume from now on that $k>1$.
Then $\G$ fixes a unique point at infinity of
the universal cover $X$ of $M$ (see e.g.~\cite{BuSc}); 
let $c(t)$ be a ray asymptotic to the point.
Since $\sec(X)$ is bounded below, the family $(X, c(t),\G)$
has a subsequence $(X,c(t_i),\G)$
that converges in the equivariant GH-topology topology to
$(X_\infty, c_\infty, \G_\infty)$.
Now $\sec(X)$ is also bounded above, the metric is regular, 
and $X$ has infinite injectivity radius, hence
the convergence $(X,c(t_i))\to (X_\infty, c_\infty)$
is in fact in $C^\infty$ topology.
Then the quotients $(X/\G,p_i)$ converge in
pointed GH-topology to $(X_\infty/\G_\infty, p_\infty)$, where
$p_i$, $p_\infty$ are the projections of $c(t_i)$, $c_\infty$,
respectively.

We now review the main results of~\cite{CFG} as they apply
to our situation; we refer to~\cite{CFG} for terminology.
Fix $\e$, $\l$ with $0<\e\ll 1\ll \l$.
By~\cite[Theorems 1.3, 1.7, Proposition 7.21]{CFG},
there are positive constants $\rho$, $\kappa$, $\nu$, $\sigma$,
depending only on $n$, $\e$, $A$ such that
for each large $i$, the manifold
$M$ carries an $N$-structure $\mathcal N_i$
and an $\mathcal N_i$-invariant $(\rho,\kappa)$-round
metric $g_{i}$ that is $\e$-close to $g$ in uniform $C^\l$-topology.
Furthermore,
there exists an orbit $O_i$ of $\mathcal N_i$ such that\newline
(i) the metric on $O_i$ induced by $g_i$ has
$\diam (O_i)\to 0$ as $i\to\infty$, \newline
(ii) $p_i$ lies in the $ \rho$-neighborhood $V_i$ of $O_i$, \newline
(iii) the normal injectivity radius of $O_i$ is $\ge \rho$,\newline
(iv) the norm of the second fundamental form of $O_i$ is $\le \nu$,
and $|\sec(O_i)|\le\s$. \newline
(v) if $\tilde V_i\to V_i$ is
the Riemannian universal cover, then $\tilde V_i$ admits a
isometric effective
action of a connected nilpotent Lie group $G_i$ that acts
transitively on the preimage $\tilde O_i$ of $O_i$
under $\tilde V_i\to V_i$, and intersects
$\pi_1(V_i)\cong\pi_1(O_i)$ in a normal
subgroup that is cocompact in $G_i$
and has index $\le \kappa$ in $\pi_1(V_i)$.

The above results are stated in~\cite{CFG} in a different form,
and their proofs are often omitted or merely sketched, so for
reader's convenience we briefly explain in the appendix how to
deduce (i)-(iv). For (v) see~\cite[pp.364--365]{CFG}.

Now we show that
the inclusion $V_i\to M$ is $\pi_1$-surjective for all large $i$.
Indeed, let $\bar V_i$ be a connected component of the preimage
of $V_i$ under the cover $X\to M$, and
as before let $\tilde V_i$, $\tilde O_i$ be the universal covers
of $V_i$, $O_i$, respectively.
Fix $\tilde q_i\in \tilde O_i$, and its projections,
$\bar q_i\in \bar V_i$ and
$q_i\in O_i\subset V_i$.
By (i)-(ii) the sequence $q_i$ subconverges to some
$q_\infty\in X_\infty/\G_\infty$, hence
for any $\g\in \G$, we have $d(\g(\bar q_i),\bar q_i)\to 0$
as $i\to \infty$. So since $\G$ is finitely generated, if $i$
is sufficiently large, then (i)-(ii)
implies that $\bar V_i$ contains the images of $\bar q_i$
under some finite generating set $S$ of $\G$.
By (iii) we see that $\bar V_i$ contains the geodesic segment
$[\bar q_i,s(\bar q_i)]$ with $s\in S$,
whose projection to $V_i\subset M$ represent
the generator of $\pi_1(M,q_i)\cong\G$ corresponding to $s$.

Hence the surjection $\pi_1(O_i)\to\G$ takes $\pi_1(O_i)\cap G_i$
onto a normal subgroup of $\G$ of index $\le\kappa$. The
intersection of all normal subgroups of $\G$ of index $\le \kappa$
is a subgroup $\G_0$ of finite index $\le n\kappa^2$ where $n=\dim
M$. (In fact, $|\G:\G_0|\le \kappa\cdot \nu_\kappa$ where
$\nu_\kappa$ is the number of normal subgroups of index $\le
\kappa$. Since $\G$ is nilpotent of $\cd(\G)<n$, it can be
generated by $<n$ elements, so there is a surjection from a rank
$n$ free group $F_n$ onto $\G$, and $\nu_\kappa$ equals to the
number of normal subgroups of $F_n$ of index $\le \kappa$, i.e.
the number of elements in $\mathrm{Hom}(F_n, \mathbb Z_\kappa)$,
which is at most $n\kappa$.)

Denote $d(\bar q_i, \g(\bar q_i))$ by $d_\g$. Below
this notation is used for different distance functions,
and each time we specify which metric we use.

Since $|\G:\G_0|<\infty$, the nilpotency degree of
$\G_0$ is $k$.
Thus there are $\g_j\in\G_0$, $j=1,\dots, k$ satisfying
\[[\g_1,[\g_2,[\g_3,[\dots[\g_{k-1},\g_k]...]=\g\ne 1.\]
Since $\G_0$ lies in the image of $\pi_1(O_i)\to\G$,
we can think of each $\g_j$ as acting on $\bar O_i\subset X$, 
where $\bar O_i$ is the preimage of $O_i$
under the cover $X\to M$.
Note that one can choose $\g_j$'s so that in the intrinsic metric
on $\bar O_i$ induced by $g_i$ we have
$d_{\g_j}\le 2n\kappa^2\cdot \diam(O_i)$.
(Indeed,  the $\G_0$-action on
$\bar O_i\subset X$ has a fundamental domain $F_i$
of diameter $\le n\kappa^2\cdot \diam(O_i)$. Then $\G_0$ is generated by
$S=\{s\in\G_0: s(\bar F_i)\cap \bar F_i\neq\emptyset\}$, and
each element of $S$ has displacement at most $2n\kappa^2\cdot \diam(O_i)$.
Then there is a nontrivial $k$-fold commutator formed by
elements of $S$, because otherwise
the identity $[a,bc]=[a,b]\cdot [b,[a,c]]\cdot [a,c]$ implies that
any $k$-fold commutator in $\G_0$ is trivial, so
its nilpotency degree is $<k$).
In particular, for the intrinsic metric induced on $O_i$ by $g_i$
the displacements of $\g_j$'s
satisfy $d_{\g_j}\to 0\ \mathrm{as}\ i\to 0$.

By (i) and (iv) we see that each $O_{i}$ with
intrinsic metric induced by $g_i$ is almost flat,
so the commutator estimate
of~\cite[Proposition~3.5 (iii), Theorem~2.4.1 (iii)]{BuKa}
for the intrinsic metric on $O_i$ induced by $g_i$ gives
\begin{equation}\label{prod}
d_{\g}\le c\prod_j d_{\g_j},
\end{equation}
where the constant $c$ depends only on $n$, $a$.

By Rauch comparison for Jacobi fields,
the normal exponential map
is bi-Lipschitz on the $\rho$-neighborhood in the normal bundle
to $O_i$, with Lipschitz constants depending on $a$, $n$, $\rho$.
Hence the nearest point projection of the $\r$-tubular
neighborhood of $\bar O_i$ onto $\bar O_i$ is $K$-Lipschitz for
$K=K(a, n, \r)$, so
any $g_i$-geodesic of length $\le 2\r$ with endpoints
on $\bar O_i$ is projected by the nearest point projection
to a curve of length $\le 2\r K$.
Since the intrinsic displacements of $\g_j$'s are $<2\r$
for all large $i$,
the estimate (\ref{prod}) holds with a different $c$,
for the distance function of the
extrinsic metric $g_{i}$, and again $c$ only depends on
$n$, $a$, $\e$, $\l$.

Finally, since the distance functions of $g$ and $g_i$
are bi-Lipschitz on $B_1(p_i))$, we get the same estimate (\ref{prod})
for the original metric $g$,
with $c$ depending on $n$, $a$, $\e$, $k$, $\l$.

For the rest of the proof we work with displacements in metric $g$.
Passing to a subsequence of $p_i$'s, we can find $j$ such that
$d_{\g_j}\ge d_{\g_l}$ for all $l$, $i$.
Taking logs we get
\[
\ln d_{\g }\le \ln C +\ln d_{\g_1}+...+\ln d_{\g_k}\le
\ln C +k\ln d_{\g_j}
\]
Since $\ln d_{\g_j}<0$ and $\lim_{i\to\infty} d_{\g_j}=0$, we deduce

\[
\limsup_{t\to\infty}\frac{\ln d_{\g}}{\ln d_{\g_j}}\ge
\limsup_{t\to\infty}\frac{\ln C}{\ln d_{\g_j}}+k=k
\]

On the other hand,
by exponential convergence of geodesic rays, for any two elements
of $\G$, and in particular for $\g,\g_j$ we get
\[
\limsup_{t\to\infty}\frac{\ln d_{\g}}{\ln d_{\g_j}}\le a+\d
\]
so $a+\d\ge k$, which completes the proof.

\begin{rmk} The weaker conclusion $a\ge k-1$ can be obtained by
the following easier argument that does not use collapsing theory.
The collapsing theory was used in the above proof to get the
commutator estimate (\ref{prod}), which is a combination of the
two independent estimates in~\cite{BuKa}, namely: \newline (a) an
upper bound on the displacement of the commutator of two elements
in terms of their displacements and rotational
parts~\cite[Corollary 2.4.2 (i)]{BuKa} that only uses bounded
curvature assumption, and \newline (b) an upper bound of the
rotational part of $\g_j$ by a constant multiple of $d_{\g_j}$
that uses almost flatness~\cite[Proposition~3.5
(i)]{BuKa}).\newline An alternative way to get (b) in our case is
via the rotation homomorphism $\phi\co \G\to O(n)$, introduced by
B.~Bowditch~\cite{Bow}, which is the holonomy of a $\G$-invariant
flat connection on $X$. A key property of $\phi$ is that
$\phi(\g)$ approximates the rotational part of any $\g\in \G$ with
error $\le d_\g$. Now since any nilpotent subgroup of $O(n)$ is
abelian, $\phi$ must have a kernel of nilpotence degree $k-1$.
Hence, there is a $(k-1)$-fold commutator in $\G$ whose entries
lie in the kernel of $\phi$, and hence their rotational parts are
bounded by their displacements. Repeating the argument at the end
of the proof of Theorem~\ref{thm: bound on pinch} for this
commutator, we get $a\ge k-1$.
\end{rmk}

\section{Infranilmanifolds are horosphere quotients}
\label{sec: optimal warp}

Let $G$ be a simply-connected nilpotent Lie group acting on itself
by left translations, and let $K$ be a compact subgroup of
$\mathrm{Aut}(G)$, so that the semidirect product $G\rtimes K$
acts on $G$ by affine transformations. The quotient of $G$ by a
discrete torsion free subgroup of $G\rtimes K$ is called an {\it
infranilmanifold}. We showed in~\cite{BK} that any pinched
negatively curved manifold with amenable fundamental group is
either the M\"obius band or product of an infranilmanifold with
$\mathbb R$, and conversely, each of these manifolds admits an
explicit warped product metrics of pinched negative curvature.

This section contains a slight improvement of
the warped product construction, that yields Theorem~\ref{optimal warp}.
Consider the product of the above $G\rtimes K$-action on $G$ with
the trivial $G\rtimes K$-action on $\mathbb R$.
For the $G\rtimes K$-action on $G\times\mathbb R$,
we prove the following.

\begin{thm}\label{thm: warp infranil}
If $G$ has nilpotence degree $k$, then for any $\e>0$,
$G\times\mathbb R$ admits a complete $G\rtimes K$-invariant
Riemannian metric of sectional curvature within $[-(k+\e)^2,-1]$.
\end{thm}
\begin{proof}
The Lie algebra $L(G)$ can
be written as
\[L(G)=L_1\supset L_2\supset\cdots\supset L_k\supset L_{k+1}=0\]
where $L_{i+1}=[L_1,L_i]$.
Note that $[L_i,L_j]\subset L_{i+j+1}$. Indeed, assume
$i\le j$ and argue by induction on $i$. The case $i=1$ is obvious
and the induction step follows from the Jacobi identity
and the induction hypothesis, because
$[L_i,L_j]=[[L_1,L_{i-1}],L_j]$ lies in
\[\mathrm{span}([[L_{i-1},L_j],L_1], [[L_1,L_j],L_{i-1}])\subset
\mathrm{span}([L_{i+j},L_1], [L_{j+1},L_{i-1}])= L_{i+j+1} \]
The group $K$ preserves each $L_i$, so we can choose a
$K$-invariant inner product $\langle\ ,\ \rangle_0$ on $L$. Let
\[F_i=\{X\in L_i\co \langle X,Y\rangle_0=0\ \mathrm{for}\ Y\in L_{i+1}\}.\]
Then $L=F_1\oplus\cdots\oplus F_k$. Define a new $K$-invariant inner
product $\langle\ ,\ \rangle_r$ on $L$ by $\langle
X,Y\rangle_r=h_i(r)^2\langle X,Y\rangle_0$ for $X,Y\in F_i$, and
$\langle X,Y\rangle_r=0$ if $X\in F_i$, $Y\in F_j$ for $i\neq j$,
where $h_i$ are some positive functions defined below. This defines
a $G\rtimes K$-invariant Riemannian metric $g_r$ on $G$.

Let $\alpha_i=i$ with $i=1,\cdots, k$ and $a=k$. Given $\rho>0$,
we define the
warping function $h_i$ to be a positive, smooth, strictly
convex, decreasing
function that is equal to $e^{-\a_i r}$ if $r\ge \rho$, and is equal to
$e^{-ar}$ if $r\le -\rho$; such a function exists since
$a\ge a_i$ for each $i$.
Thus $h_i^\prime<0<h_i^{\prime\prime}$, and the functions
$\frac{h_i^\prime}{h_i}$,
$\frac{h_i^{\prime\prime}}{h_i}$ are uniformly bounded
away from $0$ and $\infty$.

Define the warped product metric on $G\times\mathbb R$ by
$g=s^2 g_r+dr^2$, where $s>0$ is a constant; clearly $g$ is a complete
$G\rtimes K$-invariant metric.
A straightforward tedious computation
(mostly done e.g. in~\cite{BW}) yields for $g$-orthonormal vector
fields $Y_s\in F_s$ that
\begin{eqnarray*}
  & \langle R_g(Y_i,Y_j) Y_j,Y_i\rangle_g=
\frac{1}{s^2}\langle R_{g_r}(Y_i,Y_j) Y_j,Y_i\rangle_{g_r}-
\frac{h_i^\prime h_j^\prime}{h_ih_j},\\
& \langle R_g(Y_i,Y_j) Y_l,Y_m\rangle_g= \frac{1}{s^2}\langle
R_{g_r}(Y_i,Y_j) Y_l,Y_m\rangle_{g_r}
\ \ \ \mathrm{if}\ \{i,j\}\neq \{l,m\},\\
  & \langle R_g(Y_i,\frac{\partial}{\partial r})
\frac{\partial}{\partial r}), Y_i \rangle_g=
-\frac{h_i^{\prime\prime}}{h_i},\ \ \ \ \ \langle
R_g(Y_i,\frac{\partial}{\partial r})
\frac{\partial}{\partial r}), Y_j \rangle_g=0\ \ \ \mathrm{if}\ i\neq j, \\
& \langle R_g(\frac{\partial}{\partial r},Y_i) Y_j,Y_l\rangle_g=
\left(\frac{h_j^\prime}{2h_j}+\frac{h_l^\prime}{2h_l}\right)
\left(\langle [Y_j,Y_i],Y_l\rangle_g+\langle
[Y_i,Y_l],Y_j\rangle_g +\langle [Y_j,Y_l],Y_i\rangle_g\right).
\end{eqnarray*}

{\bf Correction} (added on August 28, 2010): {\it
The above formula
for $\langle R_g(\frac{\partial}{\partial r},Y_i) Y_j,Y_l\rangle_g$
is incorrect. A correction can be found in
Appendix C of~\cite{Bel-rh-warp} where it is explained
why the mistake does not affect other results
of the present paper.}

Since $[L_i,L_j]\subset L_{i+j+1}$, we have
for $Z=\sum_{i=1}^k Z_i$ and $W=\sum_{j=1}^k W_j$ with $Z_i, W_i\in F_i$
\[
|[Z,W]|_{g_r}\le \sum_{ij}|[Z_i,W_j]|_{g_r}\le\sum_{ij}
\sum_{s>i+j} h_s |[Z_i,W_j]|_{0}
\]
The above choice of $a_i$'s implies that if $r\ge \r$,
then $\sum_{s>i+j} h_s\le k h_i h_j$. Also
$|[Z_i,W_j]|_{0}\le C |Z_i|_{0}|W_j|_{0}$ where $C$ only depends on
the structure constants of $L$, so that we conclude
\[
|[Z,W]|_{g_r}\le Ck |Z_i|_{0}|W_j|_{0}
\sum_{ij} h_ih_j\le Ck^2|Z|_{g_r}|W|_{g_r}.
\]
It follows that if $r\ge \rho$, then the norm of the curvature
tensor of $g_r$ is bounded in terms of $C$, $k$~\cite[Proposition
3.18]{Cheeger-Ebin}. The same conclusions trivially hold for $r\le
-\rho$, because then $g_r$ is the rescaling of $g_0$ by a constant
$e^{-ar}>1$, and also for $r\in [-\r,\r]$ by compactness, since
$g_r$ is left-invariant and depends continuously of $r$. Hence
$\langle R_g(Y_i,Y_j) Y_l,Y_m\rangle_g\to 0$ as $s\to \infty$ if
$\{i,j\}\neq \{l,m\}$.

Also $\langle R_g(\frac{\partial}{\partial r},Y_i) Y_j,Y_l\rangle_g\to 0$
as $s\to\infty$, because
\[|\langle [Y_j,Y_i],Y_l\rangle_g|=
s^2|\langle [Y_j,Y_i],Y_l\rangle_{g_r}|\le
s^2Ck^2|Y_j|_{g_r}|Y_i|_{g_r}|Y_l|_{g_r}\le Ck^2/s,\]
where the last inequality holds since
$s|Y|_{g_r}=1$ for any $g$-unit vector $Y$.

It follows that as $s\to\infty$, then $R_g$ uniformly converges to
a tensor $\bar R$ whose nonzero components are
\[ \bar R(Y_i,Y_j, Y_j,Y_i)=
-\frac{h_i^\prime h_j^\prime}{h_ih_j}\ \ \mathrm{and}\ \ \bar
R\left(Y_i,\frac{\partial}{\partial r}, \frac{\partial}{\partial r},
Y_i\right)= -\frac{h_i^{\prime\prime}}{h_i}.\]
Thus $g$ has pinched negative curvature for all large $s$.
Finally, we show that for any $\e>0$ there exists $\r$
such that $\sec_g\in [-(k+\e)^2,-1]$. Note that
\[\frac{h_i^{\prime}}{h_i}=\ln (h_i)^\prime\qquad \mathrm{and}\qquad
\frac{h_i^{\prime\prime}}{h_i}=
\ln(h_i)^{\prime\prime}+(\ln (h_i)^\prime )^2.\]
By construction $|\ln (h_i)^\prime|\le k$. Also let $\r$ be
large enough, so that one can choose $h_i$ on $[-\r,\r]$ to satisfy
$|\ln(h_i)^{\prime\prime}|\ll \e$. Then for all sufficiently large
$s$, the sectional curvature of $g$ is within $[-(k+\e)^2,-1]$.
\end{proof}

\begin{proof}[Proof of Theorem~\ref{optimal warp}]
By~\cite{BK} if a pinched negatively curved manifold contains has
a virtually $k$-step nilpotent fundamental group, then it is
diffeomorphic to the quotient of $G\times\mathbb R$ by a discrete
torsion free subgroup of $G\rtimes K$. Thus we are done by
Theorem~\ref{thm: warp infranil}.
\end{proof}

\appendix

\section{On collapsing theory}

The purpose of this appendix is to outline the proof of the claims
(i)-(iv) made in the proof of Theorem~\ref{thm: bound on pinch}.
Some details can be found in~\cite{CFG}.

Since $g$ is regular, so is the corresponding metric $\tilde g$ on
the frame bundle $FM$. The balls $(FB_1(x), \tilde g)$ form an
$O(n)$-GH-precompact family, where $FB_1(x)$ denotes the frame
bundle over the unit ball $B_1(x)$, $x\in M$. By~\cite{Fukaya-jdg}
the closure of the family consists of regular Riemannian
manifolds. So for an arbitrary sequence $p_i\in M$, the manifolds
$(FB_1(p_i), \tilde g)$ subconverge in $O(n)$-GH-topology to a
pointed regular Riemannian manifold $(Y,y)$.

By the local version of Fukaya's
fibration theorem for some sequence $\d_i>0$ satisfying
$\d_i\to 0$ as $i\to \infty$, there exists for each large $i$
an $O(n)$-equivariant $\d_i$-almost Riemannian submersion
$FB_1(p_i)\to Y$ with nilmanifolds as fibers,
which is also an $O(n)$-$\d_i$-Hausdorff approximation.
Furthermore, each $FB_1(p_i)$ carries an $O(n)$-invariant
N-structure $\tilde{\mathcal N_i}$ whose orbits are the
nilmanifold fibers of the above submersion, and
because of the $O(n)$-invariance,
the structure descends to an N-structure $\mathcal N_i$ on $B_1(p_i)$.
By~\cite[Proposition 7.21]{CFG}
$FB_1(p_i)$ carries a metric $\tilde g_i$ that is $\e$-close
to $\tilde g$ in $C^\l$-topology, and is
both $O(n)$-invariant and $\tilde{\mathcal N_i}$-invariant.
Hence $\tilde g_i$ induces unique Riemannian submersion
metrics $\bar g_i$ on $Y$, and $g_i$ on $B_1(x)$.

To see (ii)-(iv), note that
if $l\le \l-2$, then $||\nabla^l R_{\bar g_i}||$ is bounded
independently of $i$, so the sequence $\bar g_i$ is precompact
in $C^{\l-2}$-topology.
Then by~\cite[Lemma 2.7]{PT} $\bar g_i$ is precompact
in $O(n)$-$C^{\l-2}$-topology,
i.e. after pulling back by self-diffeomorphisms
of $Y$, the metrics smoothly subconverge and share the same
isometric $O(n)$-action.
Thus there exists $\rho>0$ such that for each large $i$,
the point $y\in (Y, \bar g_i)$ lies in
a $\rho$-neighborhood of an $O(n)$-orbit that
has normal injectivity radius $\ge\rho$.
The preimage $O_i$ of the $O(n)$-orbit under the Riemannian
submersion $(FB_1(p_i),\tilde g_i)\to (Y,\bar g_i)$
satisfies (ii)-(iii). Finally, (iii) implies the second
fundamental form bound in (iv), which by Gauss formula
gives a bound on $|\sec(O_i)|$.

To see (i) note that
the $\tilde g$-diameter of any orbit
of $\tilde{\mathcal N_i}$ is $\le \d_i$, so since $\tilde g$
and $\tilde g_i$ are bi-Lipschitz, the
$\tilde g_i$-diameter of any orbit
of $\tilde{\mathcal N_i}$ tends to zero as $i\to \infty$, and
the same holds for orbits of $\mathcal N_i$ because
$FB_1(p_i)\to B_1(p_i)$ is distance nonincreasing.
Finally, the ambient diameter bound implies the intrinsic
diameter bound, because
Rauch comparison for Jacobi fields
gives bounds on bi-Lipschitz constants of
the normal exponential map of $O_i$, and in particular, the Lipschitz
constant of the nearest point projection of the $\r$-tubular
neighborhood of $O_i$ onto $O_i$ depends only on
$a$, $n$, $\r$, $\e$, $\l$.

\small
\bibliographystyle{amsalpha}
\bibliography{nilpinch}
\end{document}